\documentclass{icmart}

\contact[b.j.green@bristol.ac.uk]{School of Mathematics \\ University Walk \\ Bristol BS8 1TW \\ England}

\newtheorem{theorem}{Theorem}[section]

\newtheorem{proposition}[theorem]{Proposition}
\newtheorem{conjecture}[theorem]{Conjecture}

\newtheorem*{prin}{Principle} 

\theoremstyle{definition}
\newtheorem{definition}[theorem]{Definition}

\newtheorem{example}[subsection]{Example}
\newtheorem{question}[subsection]{Question}

\renewcommand{\leq}{\leqslant}
\renewcommand{\geq}{\geqslant}
\def\ts{\textstyle}
\def\ds{\displaystyle}

\title[Generalising the Hardy-Littlewood method for primes]{Generalising the Hardy-Littlewood method for primes}

\author[Ben Green]{Ben Green\thanks{This research was partially conducted during the period the author served as a Clay Research Fellow. He would like to express his sincere gratitude to the Clay Institute, and also to the Massachusetts Institute of Technology, where he was a Visiting Professor for the academic year 2005-06.}}

\begin{document}

\begin{abstract}
The Hardy-Littlewood method is a well-known technique in analytic number theory. Among its spectacular applications are Vinogradov's 1937 result that every sufficiently large odd number is a sum of three primes, and a related result of Chowla and Van der Corput giving an asymptotic for the number of 3-term progressions of primes, all less than $N$. This article surveys recent developments of the author and T. Tao, in which the Hardy-Littlewood method has been generalised to obtain, for example, an asymptotic for the number of 4-term arithmetic progressions of primes less than $N$.
\end{abstract}



\maketitle

\section{Introduction}
Godfrey Harold Hardy and John Edensor Littlewood wrote, in the 1920s, a famous series of papers \emph{Some problems of ``partitio numerorum''}. In these papers, whose content is elegantly surveyed by Vaughan \cite{vaughan-legacy}, they developed techniques having their genesis in work of Hardy and Ramanujan on the partition function \cite{hardy-ramanujan} to well-known questions in additive number theory such as Waring's problem and the Goldbach problem.

Papers III and V in the series, \cite{hardy-littlewood-primes1, hardy-littlewood-primes2}, were devoted to the sequence of primes. In particular it was established on the assumption of the Generalised Riemann Hypothesis that every sufficiently large odd number is the sum of three primes. In 1937 Vinogradov \cite{vinogradov-paper} made a further substantial advance by removing the need for any unproved hypothesis.

The Hardy-Littlewood-Vinogradov method may be applied to give an asymptotic count for the number of solutions in primes $p_i$ to any fixed linear equation
\[ a_1 p_1 + \dots + a_t p_t = b\]
in, say, the box $p_1,\dots,p_t \leq N$, provided that at least 3 of the $a_i$ are non-zero. This includes the three-primes result, and also the result that there are infinitely many triples of primes $p_1 < p_2 < p_3$ in arithmetic progression, due to Chowla \cite{chowla} and van der Corput \cite{vdC}.

More generally the Hardy-Littlewood method may also be used to investigate systems such as $\mathbf{A} \mathbf{p} = \mathbf{b}$, where $\mathbf{A}$ is an $s \times t$ matrix with integer entries and, potentially, $s > 1$. A natural example of such a system is given by the $(k-2) \times k$ matrix

\begin{equation}\label{ap-sys} \mathbf{A} := \begin{pmatrix} 1 & - 2 & 1 & 0 & \dots & 0 & 0 & 0 \\ 0 & 1 & -2 & 1 &\dots & 0 & 0 & 0  \\ &&&& \dots & && \\ 0 & 0 & 0 & 0 & \dots & 1 & - 2 & 1 \end{pmatrix},\end{equation} in which case a solution to $\mathbf{A}\mathbf{p} = 0$ is just a $k$-term arithmetic progression of primes.

Here, unfortunately, the Hardy-Littlewood method falters in that it generally requires $t \geq 2s + 1$. In particular it cannot be used to handle progressions of length four or longer. There are certain special systems with fewer variables which \emph{can} be handled. In this context we take the opportunity to mention a beautiful result of Balog \cite{balog}, where it is shown that for any $m$ there are distinct primes $p_1 < \dots < p_m$ such that each number $\frac{1}{2}(p_i + p_j)$ is also prime, or in other words that the system
\begin{align} \nonumber
p_1 + p_2 & = 2p_{12}\\
& \dots \nonumber \\
p_{m-1} + p_m & =  2p_{m-1,m} \label{balog}
\end{align}
has a solution in primes $p_1,\dots, p_m,p_{12},\dots, p_{m-1,m}$. There is also a result of Heath-Brown \cite{heath-brown}, in which it is established that there are infinitely many four-term progressions in which three members are prime and the fourth is either a prime or a product of two primes.

The survey of Kumchev and Tolev \cite{kumchev-survey} gives a detailed account of applications of the Hardy-Littlewood method to additive prime number theory.

The aim of this survey is to give an overview of recent work of Terence Tao and I \cite{green-tao-u3inverse,green-tao-u3mobius,green-tao-prime4aps}. Our aim, which has been partially successful, is to extend the Hardy-Littlewood method so that it is capable of handling a more-or-less arbitrary system $\mathbf{A}\mathbf{p} = \mathbf{b}$, subject to the proviso that we do not expect to be able to handle any system which secretly encodes a ``binary'' problem such as Goldbach or Twin Primes.

This is a large and somewhat technical body of work. Perhaps my main aim here is to give a guide to our work so far, pointing out ways in which the various papers fit together, and future directions we plan to take. A subsidiary aim is to focus as far as possible on key concepts, rather than on details. Of course, one would normally aim to do this in a survey article. However in our case we expect that many of these details will be substantially cleaned up in future incarnations of the theory, whilst the key concepts ought to remain more-or-less as they are.

I will say rather little about our paper \cite{green-tao-longprimeaps} establishing that there are arbitrarily long arithmetic progressions of primes. Whilst there is considerable overlap between that paper and the ideas we discuss here, those methods were somewhat ``soft'' whereas the flavour of our more recent work is distinctly ``hard''. We refer the reader to the survey of Tao in Volume I of these Proceedings, and also to the surveys \cite{green-survey,kra-survey,tao-survey-1,tao-survey-2}.

To conclude this introduction let me remark that the reader should not be under the impression that the Hardy-Littlewood method only applies to linear equations in primes, or even that this is the most popular application of the method. There has, for example, been a huge amount done on the circle of questions surrounding Waring's problem. For a survey see \cite{waring-survey}. More generally there are many spectacular results where variants of the method are used to locate integer points on quite general varieties, provided of course that there are sufficiently many variables. The reader may consult Wooley's survey \cite{wooley-survey} for more information on this.

\section{The Hardy-Littlewood heuristic}

We have stated our interest in systems of linear equations in primes. While we are still somewhat lacking in theoretical results, there are heuristics which predict what answers we should expect in more-or-less any situation.

It is natural, when working with primes, to introduce the von Mangoldt function $\Lambda : \mathbb{N} \rightarrow \mathbb{R}_{\geq 0}$, defined by
\[ \Lambda(n) := \left\{ \begin{array}{ll} \log p & \mbox{if $n = p^k$ is a prime power} \\ 0 & \mbox{otherwise}.\end{array}\right.\]
The prime powers with $k \geq 2$ make a negligible contribution to any additive expression involving $\Lambda$. Thus, for example, the prime number theorem is equivalent to the statement that
\[ \mathbb{E}_{n \leq N} \Lambda(n) = 1 + o(1).\]
Here we have used the very convenient notation of expectation from probability theory, setting $\mathbb{E}_{x \in X} := |X|^{-1} \sum_{x \in X}$ for any set $X$.

We now discuss a version of the Hardy-Littlewood heuristic for systems of linear equations in primes. Here, and for the rest of the article, we restrict attention to homogeneous systems for simplicity of exposition. 

\begin{conjecture}[Hardy-Littlewood]\label{hl-heur}
Let $\mathbf{A}$ be a fixed $s \times t$ matrix with integer entries and such that there is at least one non-zero solution to $\mathbf{A}\mathbf{x} = 0$ with $x_1,\dots,x_t \geq 0$. Then \[ \mathbb{E}_{\substack{x_1,\dots, x_t \leq N\\ \mathbf{A} \mathbf{x} = 0}} \Lambda(x_1) \dots \Lambda(x_t) = \mathfrak{S}(A)(1 + o(1))\] as $N \rightarrow \infty$, where the \emph{Singular Series} $\mathfrak{S}(A)$ is equal to a product of local factors $\prod_p \alpha_p$, where 
\[ \alpha_p := \frac{\mathbb{P}(\mathbf{x} \in \mathbb{F}_p^{\times t} | \mathbf{x} \in \mathbb{F}_p^{t}, \mathbf{A}\mathbf{x} = 0)}{\mathbb{P}(\mathbf{x} \in \mathbb{F}_p^{\times t} | \mathbf{x} \in \mathbb{F}_p^{t})}.\]
\end{conjecture}
The singular series reflects ``local obstructions'' to having solutions to $\mathbf{A} \mathbf{x} = 0$ in primes; in the simple example $\mathbf{A} = \begin{pmatrix} 1 & 9 & -27 \end{pmatrix}$, where the associated equation $p_1 + 9p_2 - 27p_3 = 0$ has no solutions, one has $\alpha_3 = 0$. A more elegant formulation of the conjecture would include a ``local obstruction at $\infty$'' $\alpha_{\infty}$, in exchange for removing the hypothesis on $\mathbf{A}$.

Chowla and van der Corput's results concerning three-term progressions of primes confirm the prediction Conjecture \ref{hl-heur} for the matrix $\mathbf{A} = \begin{pmatrix} 1 & -2 & 1\end{pmatrix}$. From this it is easy to derive an \emph{asymptotic} for the number of triples $(p_1,p_2,p_3)$,  $p_1 < p_2 < p_3 \leq N$, of primes in arithmetic progression.
\begin{theorem}[Chowla, van der Corput, \cite{chowla,vdC}] \label{prime-3aps}
The number of triples of primes $(p_1,p_2,p_3)$, $p_1 < p_2 < p_3 \leq N$, in arithmetic progression is 
\[ \mathfrak{S}_3 N^2 \log^{-3} N (1 + o(1)),\]
where
\[ \mathfrak{S}_3 := \frac{1}{2}\prod_{p\geq 3} ( 1- \frac{1}{(p-1)^2}) \approx 0.3301.\]
\end{theorem}
The singular series $\mathfrak{S}_3$ is equal to $\frac{1}{4}\mathfrak{S}(\mathbf{A})$, where $\mathbf{A} = \begin{pmatrix} 1 & -2 & 1 \end{pmatrix}$, and is also half the twin prime constant.

Certain systems $\mathbf{A}\mathbf{p} = \mathbf{b}$ should be thought of as very difficult indeed, since their understanding implies an understanding of a binary problem such as the Goldbach or twin prime problem. If $\mathbf{A}$ has the property that every non-zero vector in its row span (over $\mathbb{Q}$) has at least three non-zero entries then there is no such reason to believe that it should be fantastically hard to solve. 

\begin{definition}[Non-degenerate systems]\label{non-degenerate-def} Suppose that $s,t$ are positive integers with $t \geq s + 2$. We say that an $s \times t$ matrix $\mathbf{A}$ with integer entries is \emph{non-degenerate} if it has rank $s$, and if every non-zero vector in its row span (over $\mathbb{Q}$) has at least three non-zero entries.\end{definition}

The reader may care to check that the system \eqref{ap-sys} defining a progression of length $k$ is non-degenerate.

Our eventual goal is to prove Conjecture \ref{hl-heur} for all non-degenerate systems. This goal may be subdivided into subgoals according to the value of $s$.

\begin{conjecture}[Asymptotics for $s$ simultaneous equations]\label{mainconj}Fix a value of $s \geq 1$ and suppose that $t \geq s + 2$ and that $\mathbf{A}$ is a non-degenerate $s \times t$ matrix. Then Conjecture \ref{hl-heur} holds for the system $\mathbf{A}\mathbf{p} = 0$.
\end{conjecture}

One can also formulate an appropriate conjecture for non-homogeneous systems $\mathbf{A}\mathbf{p} = \mathbf{b}$, and one would not expect to encounter significant extra difficulties in proving it. One might also try to count prime solutions to $\mathbf{A}\mathbf{p} = 0$ in which the primes $p_i$ are subject to different constraints $p_i \leq N_i$, or perhaps are constrained to lie in a fixed arithmetic progression $p_i \equiv a_i \pmod{q_i}$. One would expect all of these extensions to be relatively straightforward.

The classical Hardy-Littlewood method can handle the case $s = 1$ of Conjecture \ref{mainconj}. Our new developments have led to a solution of the case $s = 2$. In particular we can obtain an asymptotic for the number of 4-term arithmetic progressions of primes, all less than $N$:

\begin{theorem}[G.--Tao \cite{green-tao-prime4aps}] \label{4-aps-asymptotic}
The number of quadruples of primes $(p_1,p_2,p_3,p_4)$, $p_1 < p_2 < p_3 < p_4 \leq N$, in arithmetic progression is
\[ \mathfrak{S}_4 N^2 \log^{-4} N (1 + o(1)),\]
where
\[ \mathfrak{S}_4 := \frac{3}{4} \prod_{p \geq 5} (1 - \frac{3p-1}{(p-1)^3}) \approx 0.4764.\]
\end{theorem}

\section{The Hardy-Littlewood method for primes}

\label{sec3}

The aim of this section is to describe the Hardy-Littlewood method as it would normally be applied to linear equations in primes. We will sketch the proof of Theorem \ref{prime-3aps}, the asymptotic for the number of $3$-term progressions of primes. This is equivalent to the $s=1$ case of Conjecture \ref{mainconj} for the specific matrix $\mathbf{A} = \begin{pmatrix} 1 & -2 & 1 \end{pmatrix}$. Very similar means may be used to handle the general case $s =1$ of that conjecture.

The Hardy-Littlewood method is, first and foremost, a method of harmonic analysis. The primes are studied by introducing the \emph{exponential sum} (a kind of Fourier transform)
\[ S(\theta) := \mathbb{E}_{n \leq N} \Lambda(n) e(\theta n) \] for $\theta \in \mathbb{R}/\mathbb{Z}$, where $e(\alpha) := e^{2\pi i \alpha}$. It is the appearance of the circle $\mathbb{R}/\mathbb{Z}$ here which gives the Hardy-Littlewood method its alternative name. Now it is easy to check that
\[
\mathbb{E}_{x_1,x_2,x_3 \leq N} \Lambda(x_1) \Lambda(x_2) \Lambda(x_3) 1_{x_1 - 2x_2 + x_3 = 0} = \int^1_0 S(\theta)^2 S(-2\theta) \, d\theta.
\]
whence
\begin{equation}\label{fourier-formula} 
\mathbb{E}_{\substack{x_1,x_2,x_3 \leq N \\ x_1 - 2x_2 + x_3 = 0}} \Lambda(x_1) \Lambda(x_2) \Lambda(x_3)  = (2N + O(1)) \int^1_0 S(\theta)^2 S(-2\theta) \, d\theta.
\end{equation}
The method consists of gathering information about $S(\theta)$, and then using this formula to infer an asymptotic for the left-hand side. 

The process of gathering information about $S(\theta)$ leads us to another key feature of the Hardy-Littlewood method: the realisation that one must split the set of $\theta$ into two classes, the \emph{major arcs} $\mathfrak{M}$ in which $\theta \approx a/q$ for some small $q$ and the \emph{minor arcs} $\mathfrak{m} := [0,1) \setminus \mathfrak{M}$. To see why, let us attempt some simple evaluations. First of all we note that 
\[ S(0) := \mathbb{E}_{n \leq N} \Lambda(n) = 1 + o(1),\]
this being equivalent to the prime number theorem. To evaluate $S(1/2)$, observe that almost all of the support of $\Lambda$ is on odd numbers $n$, for which $e(n/2) = -1$. Thus
\[ S(1/2) := \mathbb{E}_{n \leq N} \Lambda(n) e(n/2) = -1 + o(1).\]
The evaluation of $S(1/3)$ is a little more subtle. Most of the support of $\Lambda$ is on $n$ not divisible by 3, and for those $n$ the character $e(n/3)$ takes two values according as $n \equiv 1\!\!\pmod{3}$ or $n \equiv 2\!\!\pmod{3}$. We have
\begin{align*} S(1/3) & = e(1/3)\mathbb{E}_{n\leq N} 1_{n \equiv 1 \!\!\!\!\!\pmod 3}\Lambda(n) + e(2/3)\mathbb{E}_{n\leq N} 1_{n \equiv 2 \!\!\!\!\!\pmod 3}\Lambda(n) + o(1) \\ & = -1/2 + o(1),\end{align*} this being a consequence of the fact that the primes are asymptotically equally divided between the congruence classes $1 \!\!\pmod{3}$ and $2 \!\!\pmod{3}$.

In similar fashion one can get an estimate for $S(a/q)$ for small $q$, and indeed for $S(a/q + \eta)$ for sufficiently small $\eta$, if one uses the prime number theorem in arithmetic progressions. The set of such $\theta$ is called the \emph{major arcs} and is denoted $\mathfrak{M}$. (The notion of ``small $q$'' might be $q \leq \log^A N$, for some fixed $A$. The notion of ``small $\eta$'' might be $|\eta| \leq \log^A N/qN$. The flexibility allowed here depends on what type of prime number theorem along arithmetic progressions one is assuming. Unconditionally, the best such theorem is due to Siegel and Walfisz and it is this theorem which leads to these bounds on $q$ and $|\eta|$.) 

Suppose by contrast that $\theta \notin \mathfrak{M}$, that is to say $\theta$ is not close to $a/q$ with $q$ small. We say that $\theta \in \mathfrak{m}$, the \emph{minor arcs}. It is hard to imagine that in the sum
\begin{equation}\label{eq8} S(\sqrt{2} - 1) = \mathbb{E}_{n \leq N} \Lambda(n) e(n \sqrt{2})\end{equation}
the phases $e(n\sqrt{2})$ could conspire with $\Lambda(n)$ to prevent cancellation. It turns out that indeed there \emph{is} substantial cancellation in this sum. This was first proved by Vinogradov, and nowadays it is most readily established using an identity of Vaughan \cite{vaughan-cr}, which allows one to decompose \eqref{eq8} into three further sums which are amenable to estimation. We will discuss a variant of this method in \S \ref{vaughan-sec}. For the particular value $\theta = \sqrt{2} - 1$, and for other highly irrational values, one can obtain an estimate of the shape $|S(\theta)| \ll N^{-c}$ for some $c > 0$, which is quite remarkable since applying the best-known error term in the prime number theorem only allows one to estimate $S(0)$ with the much larger error $O(\exp(-C_{\epsilon}\log^{3/5 - \epsilon} N))$. By defining parameters suitably (that is by taking a suitable value of the constant $A$ in the precise definition of $\mathfrak{M}$), one can arrange that $S(\theta)$ is always very small indeed on the minor arcs $\mathfrak{m}$, say
\begin{equation}\label{minor-arcs-small} \sup_{\theta \in \mathfrak{m}} |S(\theta)| \ll \log^{-10} N.\end{equation} 

Recall now the formula \eqref{fourier-formula}. Splitting the integral into that over $\mathfrak{M}$ and that over $\mathfrak{m}$, we see from Parseval's identity that 
\begin{equation}\label{eq34} |\int_{\mathfrak{m}} S(\theta)^2 S(-2\theta) \, d\theta | \leq \sup_{\theta \in \mathfrak{m}}|S(\theta)| \int^1_0 |S(\theta)|^2 \, d\theta \ll \frac{\log^{-9}N}{N}.\end{equation} Thus in the effort to establish Theorem \ref{prime-3aps} the contribution from the minor arcs $\mathfrak{m}$ may essentially be ignored. The proof of that theorem is now reduced to showing that 
\[ \int_{\mathfrak{M}} S(\theta)^2 S(-2\theta)\, d\theta = (1 + o(1))\frac{1}{N} \prod_{p \geq 3} (1 - \frac{1}{(p-1)^2}).\]
Since one has asymptotic formul{\ae} for $S(\theta)$ (and $S(-2\theta)$) on $\mathfrak{M}$, this is essentially just a computation, albeit not a particularly straightforward one.

It is instructive to look for the point in the above argument where we used the fact that $\mathbf{A}$ was non-degenerate, that is to say that our problem had at least three variables. Why can we not use the same ideas to solve the twin prime or Goldbach problems? The answer lies in the bound \eqref{eq34}. In the twin prime problem we would be looking to bound
\[ |\int_{\theta \in \mathfrak{m}} |S(\theta)|^2 e(2\theta)|,\]
and the only obvious means of doing this is via an inequality of the form
\[ |\int_{\theta \in \mathfrak{m}} |S(\theta)|^2 e(2\theta)| \leq \sup_{\theta \in \mathfrak{m}} |S(\theta)|^c \int^1_0 |S(\theta)|^{2 - c} \, d\theta.\]
Now, however, Parseval's identity does not permit one to place a bound on \[ \int^1_0 |S(\theta)|^{2-c}\,d\theta.\] Indeed this whole endeavour is rather futile since heuristics predict that the minor arcs actually make a significant contribution to the asymptotic for twin primes.

An attempt to count 4-term progressions in primes via the circle method is beset by difficulties of a similar kind. 

\section{Exponential sums with M\"obius}\label{mob-sec}

The presentation in the next two sections (and in our papers) is influenced by that in the beautiful book of Iwaniec and Kowalski \cite{iwaniec-kowalski}.

In the previous section we described what is more-or-less the standard approach to solving linear equations in primes using the Hardy-Littlewood method. In \cite[Ch. 19]{iwaniec-kowalski} one may find a very elegant variant in which the M\"obius function $\mu$ is made to play a prominent r\^ole. As we saw above the behaviour of the exponential sum $S(\theta)$ was a little complicated to describe, depending as it does on how close to a rational $\theta$ is. By contrast the exponential sum
\[ M(\theta) := \mathbb{E}_{n \leq N} \mu(n) e(\theta n)\]
has a very simple behaviour, as the following result of Davenport shows.

\begin{proposition}[Davenport's Bound]\label{davenport-prop} We have the estimate \[ |M(\theta)| \ll_A \log^{-A} N\] uniformly in $\theta \in [0,1)$ for any $A > 0$.
\end{proposition}

In fact on the GRH Baker and Harman \cite{baker-harman} obtain the superior bound $|M(\theta)| \ll N^{-3/4 + \epsilon}$. By analogy with results of Salem and Zygmund \cite{salem-zygmund}  concerning random trigonometric series one might guess that the truth is that $\sup_{\theta \in [0,1)} |M(\theta)| \sim c\sqrt{\log N/N}$. This is far from known even on GRH; so far as I am aware no \emph{lower} bound of the form $\sup_{\theta \in [0,1)} |M(\theta)|\sqrt{N} \rightarrow \infty$ is known.

Although Davenport's result is easy to describe its proof has the same ingredients as used in the analysis of $S(\theta)$. One must again divide $\mathbb{R}/\mathbb{Z}$ into major and minor arcs. On the major arcs one must once more use information equivalent to a prime number theorem along arithmetic progressions, that is to say information on the zeros of $L$-functions $L(s,\chi)$ close to the line $\Re s = 1$. On the minor arcs one uses an appropriate version of Vaughan's identity. One of the attractions of working with M\"obius is that this identity takes a particularly simple form (see \cite[Ch. 13]{iwaniec-kowalski} or \cite{green-tao-u3mobius}).

We offer a rough sketch of how Proposition \ref{davenport-prop} may be used as the main ingredient in a proof of Theorem \ref{prime-3aps}, referring the reader to \cite[Ch. 19]{iwaniec-kowalski} for the details. The key point is that one has the identity
\[ \Lambda(n) = \sum_{d | n} \mu(d) \log(n/d).\]
One splits the sum over $d$ into the ranges $d \leq N^{1/10}$ and $d > N^{1/10}$ (say), obtaining a decomposition $\Lambda = \Lambda^{\sharp} + \Lambda^{\flat}$. One has
\[ S^{\flat}(\theta) :=  \mathbb{E}_{n \leq N} \Lambda^{\flat}(n) e(n\theta) = \sum_{d \leq N^{1/10}} \log d \sum_{N^{1/10} \leq k \leq N/d} \mu(k) e(\theta k d),\]
from which it follows easily using Davenport's bound that
\begin{equation}\label{lambda-flat-uniform} S^{\flat}(\theta) \ll_A \log^{-A} N\end{equation}
uniformly in $\theta \in [0,1)$.

One may then write the expression
\[ \mathbb{E}_{\substack{x_1,x_2,x_3 \leq N \\ x_1 - 2x_2 + x_3 = 0}} \Lambda(x_1)\Lambda(x_2)\Lambda(x_3)\]  as a sum of eight terms using the splitting $\Lambda = \Lambda^{\sharp} + \Lambda^{\flat}$.
The basic idea is now that the main term $\prod_{p} \alpha_p$ in Theorem \ref{prime-3aps} comes from the term with three copies of $\Lambda^{\sharp}$, whilst the other 7 terms (each of which contains at least one $\Lambda^{\flat}$) provide a negligible contribution in view of \eqref{lambda-flat-uniform} and simple variants of the formula \eqref{fourier-formula}.

We have extolled the virtues of the M\"obius function by pointing to the aesthetic qualities of Davenport's bound. A more persuasive argument for focussing on it is the following basic metaprinciple of analytic number theory:

\begin{prin}[M\"obius randomness law] The M\"obius function is highly orthogonal to any ``reasonable'' bounded function $f : \mathbb{N} \rightarrow \mathbb{C}$. That is to say
\[ \mathbb{E}_{n \leq N} \mu(n) f(n) = o(1),\]
and usually one would in fact expect 
\begin{equation}\label{strong-mob} \mathbb{E}_{n \leq N} \mu(n) f(n) \ll N^{-1/2 + \epsilon}.\end{equation}
\end{prin}

In the category ``reasonable'' in this context one would certainly include polynomials phases and other somewhat continuous objects, but one should exclude functions $f$ which are closely related to the primes ($f = \mu$ and $f = \Lambda$, for example, are clearly not orthogonal to M\"obius). 

At a finer level than is relevant to our work, the M\"obius randomness law is more reliable than other heuristics that one might formulate, for example concerning $\Lambda$. In \cite{iwaniec-luo-sarnak} it is shown that 
\[ \mathbb{E}_{n \leq N} \Lambda(n) \lambda(n) e(-2\sqrt{n}) \sim c N^{-1/4},\]
where $\lambda(n) := n^{-11/2} \tau(n)$ is a normalised version of Ramanujan's $\tau$-function. One could hardly called na\"{\i}ve for expecting square root cancellation here.

\section{Proving the M\"obius randomness law}\label{vaughan-sec}

In the last section we mentioned a principle, the \emph{M\"obius randomness law,} which is very useful as a guiding principle in analytic number theory. Unfortunately it is not possible to prove the strong version \eqref{strong-mob} of the principle in any case -- even when $f(n) \equiv 1$ it is equivalent to the Riemann hypothesis.

It is, however, possible to prove weaker estimates of the form
\begin{equation}\label{strong-asymptotic-orth} \mathbb{E}_{n \leq N} \mu(n) f(n) \ll_A \log^{-A} N,\end{equation}
for arbitrary $A > 0$, for a wide variety of functions $f$. Davenport's bound is precisely this result when $f(n) = e(\theta n)$ (and, furthermore, this result is uniform in $\theta$). Similar statements are also known for polynomial phases and for Dirichlet characters (uniformly over all characters of a fixed conductor).

Now when it comes to proving an estimate of the form \eqref{strong-asymptotic-orth}, one should think of there being two different classes of behaviour for $f$. In the first class are those $f$ which are in a vague sense multiplicative, or linear combinations of a few multiplicative functions. Then the behaviour of $\mathbb{E}_{n \leq N} \mu(n) f(n)$ can be intimately connected with the zeros of $L$-functions. One has, for example, the formula
\[ \sum_{n = 1}^{\infty} \mu(n) \chi(n)n^{-s} = \frac{1}{L(s,\chi)}\] for any fixed Dirichlet character $\chi$. By the standard contour integration technique (Perron's formula) of analytic number theory one sees that $\mathbb{E}_{n \leq N} \mu(n) \chi(n)$ is small provided that $L(s,\chi)$ does not have zeros close to $\Re s = 1$. (In fact, as reported on \cite[p. 124]{iwaniec-kowalski}, there are complications caused by possible multiple zeros of $L$, and it is better to work first with the sum $\mathbb{E}_{n \leq N} \Lambda(n) \chi(n)$ of $\chi$ over primes.)

The need to consider zeros of $L$-functions can also be felt when considering \emph{additive} characters $e(an/q)$, for relatively small $q$. Indeed any Dirichlet character to the modulus $q$ may be expressed as a linear combination of such characters. Conversely any additive character $e(an/q)$ may be written as a linear combination of Dirichlet characters to moduli dividing $q$ by using Gauss sums. By applying Siegel's theorem, which gives the best unconditional information concerning the location of zeros of $L(s,\chi)$ near to $\Re s = 1$, one obtains for any $A$ the estimate
\[ \mathbb{E}_{n \leq N} \mu(n) e(an/q) \ll_A \log^{-A} N,\]
uniformly for $q \leq \log^A N$. By partial summation the same estimate holds when $a/q$ is replaced by $\theta = a/q + \eta$ for suitably small $\eta$, that is to say for all $\theta$ which lie in the set $\mathfrak{M}$ of major arcs. 

We turn now to a completely different technique for bounding $\mathbb{E}_{n \leq N} \mu(n) f(n)$. Remarkably this is at its most effective when the previous technique fails, that is to say when $f$ is somehow \emph{far} from multiplicative.

\begin{proposition}[Type I and II sums control sums with M\"obius]\label{prop2}  Let $f : \mathbb{N} \rightarrow \mathbb{C}$ be a function with $\Vert f \Vert_{\infty} \leq 1$, and suppose that the following two estimates hold. 
\begin{enumerate}
\item \textup{(Type I sums are small)} For all $D \leq N^{2/3}$, and for all sequences $(a_d)_{d = D}^{2D}$ with $\Vert a \Vert_{l^2[D,2D)} = 1$,  we have
\begin{equation}\label{typei} |\sum_{d = D}^{2D} \sum_{1 \leq w < N/d} a_d f(wd)| \ll_A N(\log N)^{-A - 3}.\end{equation}
\item \textup{(Type II sums are small)} For all $D,W$, $N^{1/3} \leq D \leq N^{2/3}$, $N^{1/3} \leq W \leq N/D$ and all choices of complex sequences $(a_d)_{d = D}^{2D},(b_w)_{w = W}^{2W}$ with $\Vert a \Vert_{l^2[D,2D)}$ $= \Vert b \Vert_{l^2[W,2W)} = 1$, we have
\begin{equation}\label{typeii}| \sum_{d = D}^{2D} \sum_{W \leq w \leq 2W} a_d b_w f(wd)| \ll_A N(\log N)^{-A-5}.\end{equation} 
\end{enumerate}
Then \begin{equation}\label{orthog} \mathbb{E}_{n \leq N} \mu(n) f(n) \ll_A \log^{-A}N.\end{equation}
\end{proposition}

The reader may find a proof of this statement in \cite[Ch. 6]{green-tao-u3mobius}. It is proved by decomposing the M\"obius function into two parts using an identity of Vaughan \cite{vaughan-cr}. When one multiplies by $f(n)$ and sums, one of these parts leads to Type I sums and the other to Type II sums. Note that there is considerable flexibility in arranging the ranges of $D$ in which Type I and II estimates are required, but it is not important to have such flexibility in our arguments. 

The statement of Proposition \ref{prop2} may look complicated. What has been achieved, however, is the elimination of $\mu$. Strictly speaking, one actually only \emph{needs} Type I and II estimates for some rather specific choices of coefficients $a_d, b_w$ whose definition involves $\mu$. The important realisation is that it is best to forget about the precise forms of these coefficients, the general expressions \eqref{typei} and \eqref{typeii} laying bare the important underlying information required of $f$.

Note that if $f$ is close to multiplicative then there is no hope of obtaining enough cancellation in Type II sums to make use of Proposition \ref{prop2}. If $f$ is actually completely multiplicative, for example, one may take $a_d = \overline{f(d)}$ and $b_w = \overline{f(w)}$ and there is manifestly no cancellation at all in \eqref{typeii}. If this is not the case, however, then very often it \emph{is} possible to verify the bounds \eqref{typei} and \eqref{typeii}. An example of this is a linear phase $e(\theta n)$ where $\theta$ lies in the minor arcs $\mathfrak{m}$, that is to say $\theta$ is not close to $a/q$ with $q$ small. By verifying these two estimates for such $\theta$, one has from \eqref{orthog} that Davenport's bound holds when $\theta \in \mathfrak{m}$. This completes the proof of Davenport's bound, since the major arcs $\mathfrak{M}$ have already been handled using $L$-function technology.

To see how this is usually achieved in practice we refer the reader to \cite[Ch. 24]{davenport-book}. There the reader will see that a key device is the Cauchy-Schwarz inequality, which allows one to elimiate the arbitrary coefficients $a_d, b_w$.

In \cite{green-tao-u3mobius} there is also a discussion of this result. Although logically equivalent, this discussion takes a point of view which turns out to be invaluable when dealing with more complicated situations. Taking $f(n) = e(\theta n)$ in Proposition \ref{prop2}, we suppose that either \eqref{typei} or \eqref{typeii} does not hold, that is to say that either a Type I or a Type II sum is large. We then \emph{deduce} that $\theta$ must be close to a rational with small denominator, that is to say $\theta$ must be major arc. This \emph{inverse} approach to bounding sums with M\"obius means that there is no need to make an \emph{a priori} definition of what a ``major'' or ``minor'' object is. In situations to be discussed later this helps enormously.

\section{The insufficiency of harmonic analysis}

What did we mean when we stated that the Hardy-Littlewood method was a method of harmonic analysis? In \S \ref{sec3} we saw that there is a formula, \eqref{fourier-formula}, which expresses the number of 3-term progressions in a set (such as the primes) in terms of the exponential sum over that set. The following proposition is an easy consequence of a slightly generalised version of that formula:

\begin{proposition}\label{linear-1step}
Suppose that $f_1,f_2,f_3 : [N] \rightarrow [-1,1]$ are three functions and that \[ |\mathbb{E}_{\substack{x_1,x_2,x_3 \\ x_1 - 2x_2 + x_3 = 0}} f_1(x_1) f_2(x_2) f_3(x_3)| \geq \delta.\] Then for any $i = 1,2,3$ we have
\begin{equation}\label{linear-bias-0} \sup_{\theta \in [0,1)}| \mathbb{E}_{n \leq N} f_i(n) e(n\theta)| \geq (1 + o(1))\delta/2.\end{equation}
\end{proposition}

We think of this as a statement the effect that the linear exponentials $e(n\theta)$ form a \emph{characteristic system} for the linear equation $x_1 - 2x_2 + x_3 = 0$. It follows immediately from Proposition \ref{linear-1step} and Davenport's bound that M\"obius exhibits cancellation along 3-term APs, in the sense that 
\[ \mathbb{E}_{\substack{x_1,x_2,x_3 \\ x_1 - 2x_2 + x_3 = 0}} \mu(x_1)\mu(x_2)\mu(x_3) \ll_A \log^{-A} N.\]
Proposition \ref{linear-1step} is also useful for counting progressions in sets $A \subseteq [N]$, in which context one would take various of the $f_i$ to equal the \emph{balanced function} $f_A := 1_A - \alpha$ of $A$, where $\alpha := |A|/N$. It is easy to deduce from Proposition \ref{linear-1step} the following variant, which covers this situation.

\begin{proposition}\label{linear-1step-sets}
Suppose that $A \subseteq [N]$ is a set with $|A| = \alpha N$ and that \[ |\mathbb{E}_{\substack{x_1,x_2,x_3 \\ x_1 - 2x_2 + x_3 = 0}} 1_A(x_1) 1_A(x_2) 1_A(x_3) - \alpha^3 | \geq \delta.\] Write \[ f_A := 1_A - \alpha\] for the balanced function of $A$. Then we have
\begin{equation}\label{linear-bias-1} \sup_{\theta \in [0,1)}| \mathbb{E}_{n \leq N} f_A(n) e(n\theta)| \geq (1 + o(1))\delta/14.\end{equation}
\end{proposition}

If a function $f$ correlates with a linear exponential as in \eqref{linear-bias-0} or \eqref{linear-bias-1} then we sometimes say that $f$ has \emph{linear bias}.

In this section we give examples which show that the linear exponentials do not form a characteristic system for the pair of equations $x_1 - 2x_2 + x_3 = x_2 - 2x_3 + x_4 = 0$ defining a four-term progression. These examples show, in a strong sense, that the Hardy-Littlewood method in its traditional form \emph{cannot} be used to study 4-term progressions. An interesting feature of these two examples is that they were both essentially discovered by Furstenberg and Weiss \cite{furst-weiss} in the context of ergodic theory. Much of our work is paralleled in, and in fact motivated by, the work of the ergodic theory community. See the lecture by Tao in Volume 1 of these proceedings, or the elegant surveys of Kra \cite{kra-survey,kra-icm} for more discussion and references. The examples were rediscovered, in the finite setting, by Gowers \cite{gowers-4aps,gowers-longaps} in his work on Szemer\'edi's theorem.

\begin{example}[Quadratic and generalised quadratic behaviour]
Let $\alpha > 0$ be a small, fixed, real number, and define the following sets.  Let $A_1$ be defined by
\[ A_1 := \{x \in [N] : \{ x^2 \sqrt{2} \} \in [-\alpha/2,\alpha/2]\}\]
(here, $\{t\}$ denotes the fractional part of $t$, and lies in $(-1/2,1/2]$). Define also
\[ A_2 := \{x \in [N] : \{ x\sqrt{2} \{ x \sqrt{3}\}\} \in [-\alpha/2,\alpha/2]\}.\]
\end{example} Now it can be shown (not altogether straightforwardly) that $|A_1|, |A_2| \approx \alpha N$, and furthermore that 
\[ \sup_{\theta \in [0,1)}|\mathbb{E}_{n \leq N} f_{A_i} e(n\theta)| \ll N^{-c}\] for $i = 1,2$.
Thus neither of the sets $A_1,A_2$ has linear bias in a rather strong sense. If the analogue of Proposition \ref{linear-1step-sets} were true for four term progressions, then, one would expect both $A_1$ and $A_2$ to have approximately $\alpha^4N^2/6$ four-term progressions.

 The set $A_1$, however, has considerably more 4-term APs that this in view of the identity
\begin{equation}\label{constraint} x^2 - 3(x+d)^2 + 3(x+2d)^3 - (x + 3d)^2 = 0.\end{equation}
This means that if $x,x+d, x+2d \in A_1$ then 
\[ \{(x + 3d)^2 \sqrt{2}\} \in [-7\alpha/2,7\alpha/2],\] which would suggest that $x + 3d \in A_1$ with probability $\gg 1$. In fact one can show using harmonic analysis that \eqref{constraint} is the only relevant constraint in the sense that
\begin{align*} \mathbb{P}(x + 3d \in A_1 |& x, x+ d, x+2d \in A_1) \\   & \approx \mathbb{P}(y_1 - 3y_2 + 3y_3 \in [-1,1] | y_1, y_2, y_3 \in [-1,1])  = 8/27.\end{align*}
The number of 3-term progressions in $A_1$ is $\approx \alpha^3 N/4$, and so it follows that the number of 4-term progressions in $A_1$ is $\approx 2\alpha^3/27$.

The analysis of $A_2$ is rather more complicated. However one may check that if $|\{x \sqrt{3} \}| , |\{d \sqrt{3}\}|\leq 1/10$ and if $|\{y \sqrt{2} \{y \sqrt{3}\}\}| \leq \alpha/10$ for $y = x, x+d, x+2d$, then $x + 3d \in A_2$. One can show that there are $\gg \alpha^3 N^2$ choices of $x,d$ satisfying these constraints, and hence once again $A_2$ contains $\gg \alpha^3 N^2$ 4-term progressions.

\section{Generalised quadratic obstructions}

We saw in the last section that the set of linear exponentials $e(\theta n)$ is not a characteristic system for 4-term progressions. There we saw examples involving quadratics $n^2 \theta$ and generalised quadratics $n\theta_1\{n\theta_2\}$, and these must clearly be addressed by any generalisation of Propositions \ref{linear-1step} and \ref{linear-1step-sets} to 4-term APs. Somewhat remarkably, these quadratic and generalised quadratic examples are in a sense the only ones.

\begin{proposition}\label{linear-2step}
Suppose that $f_1,f_2,f_3,f_4 : [N] \rightarrow [-1,1]$ are four functions and that \begin{equation}\label{eq299} |\mathbb{E}_{\substack{x_1,x_2,x_3,x_4 \\ x_1 - 2x_2 + x_3 = 0 \\ x_2 - 2x_3 + x_4 = 0}} f_1(x_1) f_2(x_2) f_3(x_3) f_4(x_4)| \geq \delta.\end{equation} Then for any $i = 1,2,3,4$ there is 
a generalised quadratic polynomial
\begin{equation}\label{phi-form} \phi(n) = \sum_{r,s \leq C_1(\delta)} \beta_{rs} \{ \theta_r n\} \{ \theta_s n\} + \gamma_r \{ \theta_r n\},\end{equation}
where $\beta_{rs},\gamma_r, \theta_r \in \mathbb{R}$, such that 
\[ |\mathbb{E}_{n \leq N} f_i(n) e(\phi(n))| \geq c_2(\delta).\]
We can take $C_1(\delta) \sim \exp(\delta^{-C})$ and $c_2(\delta) \sim \exp(-\delta^{-C})$.
\end{proposition}

Note that
\[ \theta n^2 = 100\theta N^2\{\frac{n}{10N}\}^2\]
and
\[ \theta_1 n \{\theta_2 n\} = 10\theta_1 N \{ \frac{n}{10N} \} \{ \theta_2 n\}\]
for $n \leq N$, and so the phases which can be written in the form \eqref{phi-form} do include all those which were discovered to be relevant in the preceding section.

The proof of Proposition \ref{linear-2step} is given in \cite{green-tao-u3inverse}. It builds on earlier work of Gowers \cite{gowers-4aps,gowers-longaps}. In \cite{green-tao-u3inverse} (see also \cite{green-tao-u3mobius}) several results of a related nature are given, in which other characteristic systems for the equation $x_1 - 2x_2 + x_3 = x_2 -2x_3 + x_4 = 0$ are given. These systems all have a ``quadratic'' flavour. We will discuss the family of \emph{$2$-step nilsequences}, which is perhaps the most conceptually appealing, in \S \ref{sec9}. In \S \ref{sec10} we will mention the family of \emph{local quadratics}, which are useful for computations involving the M\"obius function. The only real merit of the generalised quadratic phases $e(\phi(n))$ discussed above is that they are easy to describe from first principles.

\section{The Gowers norms and inverse theorems}

The proof of Proposition \ref{linear-2step} is long and complicated: there does not seem to be anything so simple as Formula \eqref{fourier-formula} in the world of 4-term progressions. Very roughly speaking one assumes that \eqref{eq299} holds, and then one proceeds to place more and more structure on each function $f_i$ until eventually one establishes that $f_i$ correlates with a generalised quadratic phase. There is a finite field setting for this argument, and we would recommend that the interested reader read this first: it may be found in \cite[Ch. 5]{green-tao-u3inverse}. The ICM lecture of Gowers \cite{gowers-icm} is a fine introduction to the ideas in his paper \cite{gowers-4aps}, which is the foundation of our work. 

There is only one part of the existing theory which we feel sure will play some r\^ole in future incarnations of these methods. This is the first step in the long series of deductions from \eqref{eq299}, in which one shows that each $f_i$ has large \emph{Gowers norm}. For the purposes of this exposition\footnote{In practice we do all our work the group $\mathbb{Z}/N'\mathbb{Z}$ for some prime  $N' \geq N$ with $N' \approx M(\mathbf{A})N$, where $M(\mathbf{A})$ is some constant depending on the system of equations $\mathbf{A}\mathbf{x} = 0$ one is interested in. One advantage of this is that the number of solutions to $\mathbf{A}\mathbf{x} = 0$ in $\mathbb{Z}/N'\mathbb{Z}$ is much easier to count than the number of solutions in $[N]$. The Gowers norms defined here differ from the Gowers norms in those settings by constant factors, so for expository purposes they may be thought of as the same. In the group setting the constant $c_{\mathbf{A}}$ in Proposition \ref{gvn-1} is simply 1.} we define the Gowers
$U^2$-norm $\Vert f \Vert_{U^2}$ of a function $f : [N] \rightarrow [-1,1]$ by
\[ \Vert f \Vert_{U^2}^4 := \mathbb{E}_{\substack{x_{00}, x_{01}, x_{10}, x_{11} \leq N \\ x_{00} + x_{11} = x_{01} + x_{10}}} f(x_{00}) f(x_{01}) f(x_{10}) f(x_{11}),\]
which is a sort of average of $f$ over two dimensional parallelograms. The $U^k$ norm, $k \geq 3$, is an average of $f$ over $k$-dimensional parallelepipeds. Written down formally it looks much more complicated than it is:
\[ \Vert f \Vert_{U^k}^{2^k} := \mathbb{E}_{\substack{x_{0,\dots,0}, \dots, x_{1,\dots,1} \\ x_{\omega^{(1)}} + x_{\omega^{(2)}} = x_{\omega^{(3)}} + x_{\omega^{(4)}}}} f(x_{0,\dots,0}) \dots f(x_{1,\dots,1}),\]
where there are $2^k$ variables $x_{\omega}$, $\omega = (\omega_1,\dots,\omega_k) \in \{0,1\}^k$, and the constraints range over all quadruples $(\omega^{(1)}, \omega^{(2)}, \omega^{(3)}, \omega^{(4)}) \in (\{0,1\}^k)^4$ with $\omega^{(1)} + \omega^{(2)} = \omega^{(3)} + \omega^{(4)}$.

The Gowers $U^k$ norm governs the behaviour of any non-degenerate system $\mathbf{A}\mathbf{x} = 0$ in which $\mathbf{A}$ has $(k-1)$ rows.

\begin{proposition}[Generalised von Neumann theorem]\label{gvn-1}
Suppose that $\mathbf{A}$ is a non-degenerate $s \times t$ matrix with integer entries. Suppose that $f_1,\dots,f_t : [N] \rightarrow [-1,1]$ are functions and that 
\[ |\mathbb{E}_{\substack{x_1,\dots,x_t \\ \mathbf{A}\mathbf{x} = 0}} f_1(x_1) \dots f_t(x_t)| \geq \delta.\]
Then for each $i = 1,\dots,t$ we have
\[ \Vert f_i \Vert_{U^{s+1}} \geq c_{\mathbf{A}}\delta.\]
\end{proposition}
The proof involves $s+1$ applications of the Cauchy-Schwarz inequality. In this generality, the result was obtained in \cite{green-tao-u3mobius}, though the proof technique is the same as in \cite{gowers-longaps}. There are results in ergodic theory of the same general type, in which ``non-conventional ergodic averages'' are bounded using seminorms which are analogous to the $U^k$-norms: see \cite{host-kra}.

Taking $s = k-2$ and $\mathbf{A}$ as in \eqref{ap-sys}, we see that in particular the Gowers $U^{k-1}$-norm ``controls'' $k$-term progressions. The Gowers norms are, of course, themselves defined by a system of linear equations, and so they must be studied as part of a generalised Hardy-Littlewood method with as broad a scope as we would like. The Generalised von Neumann Theorem may be regarded as a statement to the effect that in a sense they represent the \emph{only} systems of equations that need to be studied.

The Gowers norms do not feature in the classical Hardy-Littlewood method. It is, however, possible to prove a somewhat weaker version of Proposition \ref{linear-1step} by combining the case $k = 3$ of Proposition \ref{gvn-1} with the following \emph{inverse theorem}:

\begin{proposition}[Inverse theorem for $U^2$]\label{inverseu2}
Suppose that $N$ is large and that $f : [N] \rightarrow [-1,1]$ is a function with $\Vert f \Vert_{U^2}\geq \delta$. Then we have
\[ \sup_{\theta \in [0,1)} |\mathbb{E}_{n \leq N} f(n) e(n\theta)| \geq 2\delta^2.\]
\end{proposition}
To prove this we note the formula
\[ \mathbb{E}_{x_{00},x_{01},x_{10},x_{11}} f(x_{00})f(x_{01})f(x_{10})f(x_{11}) 1_{x_{00} + x_{11} = x_{01} + x_{10}} = \int^1_0 |\widehat{f}(\theta)|^4\,d\theta,\]
where $\widehat{f}(\theta) := \mathbb{E}_{n \leq N} f(n) e(n\theta)$. This implies that
\[ \Vert f \Vert_{U^2}^4 = (3N + O(1)) \Vert \widehat{f} \Vert_4^4.\] In view of the fact that $\Vert \widehat{f} \Vert_2^2 \leq 1/N$, this and the assumption that $\Vert f \Vert_{U^2} \geq \delta$ imply that 
\[ \Vert \widehat{f} \Vert_{\infty}^2  \geq (3 + o(1)) \delta^4,\]
which implies the result.

This argument should be compared to the argument in \eqref{eq34}, to which it corresponds rather closely.

To deduce Proposition \ref{linear-1step} by passing through Proposition \ref{inverseu2} is rather perverse, since the derivation is longer than the one that proceeds via an analogue of \eqref{fourier-formula} and it leads to worse dependencies. With our current technology, however, this is the only method which is amenable to generalisation.

Similarly, one may deduce Proposition \ref{linear-2step} from Proposition \ref{gvn-1} and the following result.

\begin{proposition}[Inverse theorem for the $U^3$-norm]\label{inverseu3}
Suppose that $f : [N] \rightarrow \mathbb{R}$ is a function for which $\Vert f \Vert_{\infty} \leq 1$ and $\Vert f \Vert_{U^3} \geq \delta$. Then there is a generalised quadratic phase
\begin{equation} \phi(n) = \sum_{r,s \leq C_1(\delta)} \beta_{rs} \{ \theta_r n\} \{ \theta_s n\} + \gamma_r \{ \theta_r n\},\end{equation}
where $\beta_{rs},\gamma_r, \theta_r \in \mathbb{R}$, such that 
\[ |\mathbb{E}_{n \leq N} f(n) e(\phi(n))| \geq c_2(\delta).\]
We can take $C_1(\delta) \sim \exp(\delta^{-C})$ and $c_2(\delta) \sim \exp(-\delta^{-C})$.
\end{proposition}

This result (and variations of it involving other ``quadratic families'') is in fact the main theorem in \cite{green-tao-u3inverse}.

As we mentioned, one may find a series of seminorms which are analogous to the Gowers norms in the ergodic-theoretic work of Host and Kra \cite{host-kra}. There are no such seminorms in the related work of Ziegler \cite{ziegler}, however, and this suggests that (as in the classical case) the Gowers norms may not be completely fundamental to a generalised Hardy-Littlewood method.

\section{Nilsequences}\label{sec9}

In the previous section we introduced the Gowers $U^k$-norms, and stated inverse theorems for the $U^2$- and $U^3$- norms. These inverse theorems provide lists of rather algebraic functions which are \emph{characteristic} for a given system of equations $\mathbf{A}\mathbf{x} = 0$. Roughly speaking, the linear phases $e(\theta n)$ are characteristic for single linear equations in which $\mathbf{A}$ is a $1 \times t$ matrix. Generalised quadratic phases $e(\phi(n))$ are characteristic for pairs of linear equations in which $\mathbf{A}$ is a non-degenerate $2 \times t$ matrix.

These two results leave open the question of whether there is a similar list of functions which is characteristic for the $U^{k}$-norm, $k \geq 4$ and hence, by the Generalised von Neumann Theorem, for non-degenerate systems defined by an $s \times t$ matrix with $s \geq 3$. The form of Propositions \ref{inverseu2} and \ref{inverseu3} does not suggest a particularly natural form for such a result, however, and indeed Proposition \ref{inverseu2} is already rather unnatural-looking.

To make more natural statements, we introduce a class of functions called \emph{nilsequences}.

\begin{definition}
Let $G$ be a connected, simply connected, $k$-step nilpotent Lie group. That is, the central series $G_0 := G$, $G_{i+1} = [G, G_i]$ terminates with $G_k = \{e\}$. Let $\Gamma \subseteq G$ be a discrete, cocompact subgroup. The quotient $G/\Gamma$ is then called a $k$-step nilmanifold. The group $G$ acts on $G/\Gamma$ via the map $T_g(x \Gamma) = xg\Gamma$. If $F : G/\Gamma \rightarrow \mathbb{C}$ is a bounded, Lipschitz function and $x \in G/\Gamma$ then we refer to the sequence $(F(T_g^n \cdot x))_{n \in \mathbb{N}}$ as a $k$-step nilsequence. 
\end{definition}

By analogy with the results of Host and Kra \cite{host-kra} in ergodic theory, we expect the collection of $(k-1)$-step nilsequences to be characteristic for the $U^k$-norm. The following conjecture is one of the guiding principles of the generalised Hardy-Littlewood method.

\begin{conjecture}[Inverse conjecture for $U^k$-norms]\label{conj}
Suppose that $k \geq 2$ and that $f : [N] \rightarrow [-1,1]$ has $\Vert f \Vert_{U^k} \geq \delta$. Then there is a $(k-1)$-step nilmanifold $G/\Gamma$ with dimension at most $C_{1,k}(\delta)$, together with a function $F : G/\Gamma \rightarrow \mathbb{C}$ with $\Vert F \Vert_{\infty} \leq 1$ and Lipschitz constant at most $C_{2,k}(\delta)$ and elements $g \in G$, $x \in G/\Gamma$ such that
\begin{equation}\label{nil-cor} |\mathbb{E}_{n \leq N} f(n) F(T_g^n \cdot x)| \geq c_{3,k}(\delta).\end{equation}
\end{conjecture}

We can at least be sure that Conjecture \ref{conj} is no more complicated than necessary, since in \cite[Ch. 12]{green-tao-u3inverse} we showed that if a bounded function $f$ correlates with a $(k-1)$-step nilsequence as in \eqref{nil-cor} then $f$ \emph{does} have large Gowers $U^k$-norm. This, incidentally, is another reason to believe that the Gowers norms play a fundamental r\^ole in the theory. It is not the case that correlation of a function $f$ with a $(k-1)$-step nilsequence prohibits $f$ from enjoying cancellation along $k$-term arithmetic progressions, for example. In the case $k=3$ an example of this phenomenon is given by the function $f$ which equals $\alpha$ for $1 \leq n \leq N/3$ and $-1$ for $N/3 < n \leq N$, where $\alpha$ is the root between $1$ and $2$ of $\alpha^3 - \alpha^2 + 3\alpha - 4 = 0$. This $f$ correlates with the constant nilsequence $1$ yet exhibits cancellation along 3-term progressions, as the reader may care to check.

Conjecture \ref{conj} seems, at first sight, to be completely unrelated to Propositions \ref{inverseu2} and \ref{inverseu3}. However after a moment's thought one realises that a linear phase $e(\theta n)$ can be regarded as a $1$-step nilsequence in which $G = \mathbb{R}$ , $\Gamma = \mathbb{Z}$, $g = \theta$ and $x = 0$. Thus Proposition \ref{inverseu2} immediately implies the case $k = 2$ of Conjecture \ref{conj}. 

The case $k = 3$ is proved in \cite{green-tao-u3inverse}. One first proves Proposition \ref{inverseu3}, and then one shows how any generalised quadratic phase $e(\phi(n))$ may be approximated by a $2$-step nilsequence. Let us discuss a simple example, the \emph{Heisenberg nilmanifold}, to convince the reader that $2$-step nilsequences can give rise to ``generalised quadratic'' behaviour.

\begin{example}[The Heisenberg nilmanifold]\label{nil3}  Consider
\[ G := \begin{pmatrix}
1 & \mathbb{R} & \mathbb{R}\\
0 & 1  & \mathbb{R}\\
0 & 0  & 1
\end{pmatrix}; \quad
\Gamma := \begin{pmatrix}
1 & \mathbb{Z} & \mathbb{Z}\\
0 & 1  & \mathbb{Z}\\
0 & 0  & 1
\end{pmatrix}. \]
Then $G/\Gamma$ is a 2-step nilmanifold.  By using the identification
\[ (x,y,z) \equiv \begin{pmatrix}
1 & x  & y\\
0 & 1  & z\\
0 & 0  & 1
\end{pmatrix} \Gamma,\]
we can identify $G/\Gamma$ (as a set) with $\mathbb{R}^3$, quotiented out by the equivalence relations
\[ (x,y,z) \sim (x+a, y+b+cx, z+c) \hbox{ for all } a,b,c \in \mathbb{Z}.\]
This can in turn be coordinatised by the cylinder $(\mathbb{R}/\mathbb{Z})^2 \times [-1/2,1/2]$ with the identification
$(x,y,-1/2) \sim (x,x+y,1/2)$. Let $F : G/\Gamma \rightarrow \mathbb{C}$ be a function. We may lift this to a function $\widetilde{F} : G \rightarrow \mathbb{C}$, defined by $\widetilde{F}(g) := F(g\Gamma)$. In coordinates, this lift takes the form
\[ \widetilde{F}(x,y,z) = F(x\!\!\!\!\!\pmod{1}, y - [z] x\!\!\!\!\!\pmod{1}, \{z \})\]
where $[z] = z - \{z\}$ is the nearest integer to $x$.
Let
\[ g := \begin{pmatrix}
1 & \alpha  & \beta \\
0 & 1  & \gamma \\
0 & 0  & 1
\end{pmatrix} \]
be an element of $G$. Then the shift $T_g : G \rightarrow G$ is given by
\[ T_g(x,y,z) = (x + \alpha, y + \beta + \gamma x, z + \gamma).\]
A short induction confirms, for example, that 
\[ T_g^n(0,0,0) = (n\alpha, n\beta  + \ts\frac{1}{2}\ds n(n+1)\alpha\gamma , n\gamma).\]
 Therefore if $F : G/\Gamma \rightarrow G/\Gamma$ is any Lipschitz function, written as a function $F: (\mathbb{R}/\mathbb{Z})^2 \times [-1/2,1/2] \to \mathbb{C}$ 
with $F(-1/2,y,z) = F(1/2,y+z,z)$, then we have
\begin{align*}  F(T^n_g &(0,0,0)) \\ & = F( n\alpha \!\!\!\!\!\pmod{1}, n\beta  + \ts\frac{1}{2}\ds n(n+1) \alpha\gamma - [n\gamma] n\alpha\!\!\!\!\!\pmod{1}, \{n \gamma\}).\end{align*}

The term $[n\gamma]n\alpha$ which appears here certainly exhibits a sort of generalised quadratic behaviour. For a complete description of how an arbitrary generalised quadratic phase $e(\phi(n))$ can be approximated by a two-step nilsequence, we refer the reader to \cite[Ch. 12]{green-tao-u3inverse}.
\end{example}
Let us conclude this section by stating, for the reader's convenience, a result/conjecture which summarises much of our discussion so far in one place.

\begin{theorem}[G.--Tao \cite{green-tao-u3inverse}]\label{th-so-far} We have the following two statements.
\begin{enumerate}
\item[\textup{(i)}] \textup{(Generalised von Neumann)} Suppose that $s,t$ are positive integers with $s + 2\leq t$. Suppose that $\mathbf{A}$ is non-degenerate $s \times t$ matrix with integer entries. Suppose that $f_1,\dots,f_t : [N] \rightarrow [-1,1]$ are functions and that 
\begin{equation}\label{correlation} |\mathbb{E}_{\substack{x_1,\dots,x_t \\ \mathbf{A}\mathbf{x} = 0}} f_1(x_1) \dots f_t(x_t)| \geq \delta.\end{equation}
Then for each $i = 1,\dots,t$ we have
\[ \Vert f_i \Vert_{U^{s+1}} \geq c_{\mathbf{A}}\delta.\]
\item[\textup{(ii)}] \textup{(Gowers inverse result: proved for $k = 2,3$, conjectural for $k \geq 4$)} Suppose that $f : [N] \rightarrow [-1,1]$ has $\Vert f \Vert_{U^k} \geq \delta$. Then there is a $(k-1)$-step nilmanifold $G/\Gamma$ with dimension at most $C_{1,k}(\delta)$, together with a function $F : G/\Gamma \rightarrow \mathbb{C}$ with $\Vert F \Vert_{\infty} \leq 1$ and Lipschitz constant at most $C_{2,k}(\delta)$ and elements $g \in G$, $x \in G/\Gamma$ such that
\begin{equation}\label{nil-cor-2} |\mathbb{E}_{n \leq N} f(n) F(T_g^n \cdot x)| \geq c_{3,k}(\delta).\end{equation}
\end{enumerate}
\end{theorem}

In particular when $s = 1$ or $2$ and \eqref{correlation} holds for some $\mathbf{A}$ and some $\delta$ then for each $i = 1,\dots,t$ there is a $2$-step nilsequence $(F(T_g^n \cdot x))_{n \in \mathbb{N}}$ such that
\begin{equation}\label{nil-cor-2a} |\mathbb{E}_{n \leq N} f_i(n) F(T_g^n \cdot x)| \geq c_{\mathbf{A}}(\delta).\end{equation}

\section{Working with the primes}

Let us suppose that we wish to count four-term progressions in the primes. One might try to apply Theorem \ref{th-so-far} with the functions $f_i$ equal to the balanced function of $A$, the set of primes $p \leq N$, and then hope to rule out a correlation such as \eqref{nil-cor-2} for some $\delta = o(\alpha^t)$ (here, of course, $\alpha \approx \log^{-1}N$ by the prime number theorem). This would then lead to an asymptotic using various instances of \eqref{correlation} together with the triangle inequality. 

There are two reasons why this is a hopeless strategy. First of all, the primes \emph{do} correlate with nilsequences. In fact since all primes other than 2 are odd it is easy to see that 
\[ \mathbb{E}_{n \leq N} f_A(n) e(n/2) \approx -\alpha.\]
There is a way to circumvent this problem, which we call the $W$-trick. The idea is that if $W = 2 \times 3 \times \dots \times w(N)$ is the product of the first several primes, then for any $b$ coprime to $W$ the set
\[ A_b :=  \{ n \leq N : Wn + b \;\; \mbox{is prime}\}\]
does not exhibit significant bias in progressions with common difference $q \leq w(N)$. One can then count 4-term progressions in the primes by counting 4-term progressions in $A_{b_1} \times \dots A_{b_4}$ for each quadruple $(b_1,\dots,b_4) \in (\mathbb{Z}/W\mathbb{Z})^{\times 4}$ in arithmetic progression and adding.

We refer to any set $A_b$ as a set of ``$W$-tricked primes''. In practice one is only free to take $w(N) \sim \log \log N$, since one must be able to understand the distribution of primes in progressions with common difference $W$ (note that even on GRH one could only take $w(N) \sim c\log N$). Even assuming we could obtain optimal results concerning the correlation of the $W$-tricked primes with $2$-step nilsequences, this information will be very weak indeed.

This highlights a more serious problem with the suggested strategy. Suppose that $A \subseteq [N]$ is a set of density $\alpha$ for which there is no obvious reason why $A$ should have an unexpectedly large or small number of 4-term APs, that is to say for which we might hope to prove that 
\begin{equation}\label{eq476} \mathbb{E}_{\substack{x_1,x_2,x_3,x_4 \\ x_1 - 2x_2 + x_3 = 0 \\ x_2 - 2x_3 + x_4 = 0}} 1_A(x_1) 1_A(x_2) 1_A(x_3) 1_A(x_4) \approx \alpha^4.\end{equation}
For example, $A$ might be the $W$-tricked primes less than $N$, in which case $\alpha \sim \frac{W}{\phi(W)}\log^{-1}N$.

We might prove \eqref{eq476} by writing $1_A = \alpha + f_A$, expanding as the sum of sixteen terms, and showing that fifteen of these are $o(\alpha^4)$ by appealing to Theorem \ref{th-so-far}, and ruling out a correlation with a $2$-step nilsequence as in \eqref{nil-cor-2a}. Unfortunately we will be operating with $\delta = o(\alpha^4) \ll \log^{-4 + \epsilon}N$, and the dependence of $c_{\mathbf{A}}(\delta)$ on $\delta$ is very weak, 
being of the form $\exp(-\delta^{-C})$. Thus we are asking to rule out the possiblility that
\[ |\mathbb{E}_{n \leq N} f_A(n) F(T_g^n \cdot x) | \gg \exp(-\log^C N)\]
for some potentially rather large $C$. This is a problem, since one would never expect more than square root cancellation in any such expression. In fact for the $W$-tricked primes one only has a small amount (depending on $w(N)$) of potential cancellation to work with and to all intents and purposes one should not bank on having available any estimate stronger than
\[ \mathbb{E}_{n \leq N} f_A(n) F(T_g^n \cdot x) = o(1).\]

What one really needs is a version of Proposition \ref{eq299} which applies to functions which need not be bounded by 1. Then one could hope to work with the von Mangoldt function $\Lambda$ instead of the far less natural characteristic function $1_A$, or more accurately with $W$-tricked variants of the von Mangoldt function such as
\[ \Lambda_{b,W}(n) := \frac{\phi(W)}{W} \Lambda(Wn + b).\]
Such a result is the main result of our forthcoming paper \cite{green-tao-prime4aps}. It would take us too far afield to say anything concerning its proof, other than that it uses one of the key tools from our paper \cite{green-tao-longprimeaps} on long progressions of primes, the ``ergodic transference'' technology of \cite[Chs. 6,7,8]{green-tao-longprimeaps}.

\begin{proposition}[Transference principle, \cite{green-tao-prime4aps}]
Suppose that $\nu : [N] \rightarrow \mathbb{R}^{+}$ is a \emph{pseudorandom measure}. Then
\begin{enumerate}
\item[\textup{(i)}] The generalised von Neumann theorem, \textup{Theorem \ref{th-so-far} (i)}, continues to hold for functions $f_1,\dots,f_t : [N] \rightarrow \mathbb{R}^+$ such that $|f_i(x)| \leq 1 + \nu(x)$ pointwise \textup{(}the value of $c_{\mathbf{A}}$ may need to be reduced slightly\textup{)}.
\item[\textup{(ii)}] If the Gowers inverse conjecture, \textup{Theorem \ref{th-so-far} (ii)}, holds for a given value of $k$ then it continues to hold for a function $f$ such that $|f(x)| \leq 1 + \nu(x)$ pointwise. In particular such an extension of the Gowers inverse conjecture is true when $k = 2,3$.
\end{enumerate}
\end{proposition}

The reader may consult \cite[Ch. 3]{green-tao-longprimeaps} for a definition of the term \emph{pseudorandom measure} and a discussion concerning it. For the purposes of this article the reader can merely accept that there is such a notion, and furthermore that one may construct a pseudorandom measure $\nu : [N] \rightarrow \mathbb{R}^+$ such that $\nu + 1$ dominates any fixed $W$-tricked von Mangoldt function $\Lambda_{W,b}$. The construction of $\nu$ comes from sieve theoretic ideas originating with Selberg. The confirmation that $\nu$ is pseudorandom is essentially due, in a very different context, to Goldston and Y{\i}ld{\i}r{\i}m \cite{gy}.

Applying these two results, one may see that the Hardy-Littlewood conjecture \ref{hl-heur} for a given non-degenerate $s \times t$ matrix $\mathbf{A}$ is a consequence of the Gowers inverse conjecture in the case $k = s+1$ together with a bound of the form
\begin{equation}\label{lam-pseudo-random}
\mathbb{E}_{n \leq N} (\Lambda_{b,W} - 1) F(T_g^n \cdot x) = o_{G/\Gamma,F}(1)
\end{equation}
for every $s$-step nilsequence $(F(T_g^n \cdot x))_{n \in \mathbb{N}}$.

By effecting a decomposition of $\Lambda_{b,W}$ as $\Lambda_{b,W}^{\sharp} + \Lambda_{b,W}^{\flat}$ rather like that in \S \ref{mob-sec}, the proof of this statement may be further reduced to a similar result for the M\"obius function:

\begin{conjecture}[M\"obius and nilsequences]\label{mob-conj}
For all $A > 0$. We have the bound
\[ \mathbb{E}_{n \leq N} \mu(n) F(T_g^n \cdot x) \ll_{A,G/\Gamma,F}\log^{-A} N\]
for every $k$-step nilsequence $(F(T_g^n \cdot x))_{n \in \mathbb{N}}$.  
\end{conjecture}

Note that we require more cancellation (a power of a logarithm) here than in \eqref{lam-pseudo-random}. This is because in passing from $\mu$ to $\Lambda^{\flat}_{b,W}$ one loses a logarithm in performing partial summation as in the derivation of \eqref{lambda-flat-uniform}. The method we have in mind to prove Conjecture \ref{mob-conj}, however, is likely to give this strong cancellation at no extra cost.

Conjecture \ref{mob-conj} posits a rather vast generalisation of Davenport's bound. The conjecture is, of course, highly plausible in view of the M\"obius randomness law. 

Let us remark that the derivation of \eqref{lam-pseudo-random} from Conjecture \ref{mob-conj} is not at all immediate, since one must also handle the contribution from $\Lambda^{\sharp}_{b,W}$. To do this one uses methods of classical analytic number theory rather similar to those of Goldston and Y{\i}ld{\i}r{\i}m \cite{gy}.

\section{M\"obius and nilsequences}\label{sec10}

The main result of \cite{green-tao-u3mobius} is a proof of Conjecture \ref{mob-conj} in the case $k = 2$. This leads, by the reasoning outlined in the previous section, to a proof of Conjecture \ref{mainconj} in the case $s = 2$.

We remarked that the classical Hardy-Littlewood method was a technique of harmonic analysis. We also highlighted the idea of dividing into major and minor arcs. We have said much on the subject of generalising the underlying harmonic analysis, but as yet there has been nothing said about a suitable extension of major and minor arcs. In this section we describe such an extension by making some remarks concerning the proof of the case $k = 2$ of Conjecture \ref{mob-conj}.

In \S \ref{vaughan-sec} we discussed how bounds on Type I and II sums may be used to show that a given function $f$ does not correlate with M\"obius. Recalling our ``inverse'' strategy for proving Davenport's bound, one might be tempted to go straight into Proposition \ref{prop2} with $f(n) = F(T_g^n  \cdot x)$, a $2$-step nilsequence, posit largeness of either a Type I or a Type II sum, and then use this to say that the nilsequence is somehow ``major arc''. One might then hope to handle the major nilsequences by some other method, perhaps the theory of $L$-functions.

Such an attempt is a little too simplistic, for the following reason. Returning to the $1$-step case, note that the sum of two $1$-step nilsequences is also a $1$-step nilsequence (on the product nilmanifold $G_1/\Gamma_1 \times G_2/\Gamma_2$). In particular, the function $f(n) = e(n/5) + e(n\sqrt{2})$ is a 1-step nilsequence. We know, however, that to handle correlation of M\"obius with $e(n/5)$ we need to know something about $L$-functions, whereas we do not have an $L$-function method of handling $e(n\sqrt{2})$. This suggests that some sort of preliminary decomposition of the function $f$ is in order, and such a suggestion turns out to be correct.

In the $2$-step case, a nilsequence $F(T_g^n \cdot x)$ can be decomposed into \emph{local quadratics}. These are objects of the form
\begin{equation}\label{local-quad-def}
f(n) := 1_{B_N}(n) e(\phi(n)),
\end{equation}
where $B_N$ is a set of the form
\[ B_N := \{n : N/2 \leq n < N : F_1(n) \neq 0\}\]
for some $1$-step nilsequence $F_1$ depending on $F, G/\Gamma, g$ and $x$, and $\phi : B_N \rightarrow \mathbb{R}/\mathbb{Z}$ is \emph{locally quadratic}. This means that one may unambiguously define the second derivative $\phi''(h_1,h_2)$ to equal 
\[ \phi(x + h_1 + h_2) - \phi(x+h_1) - \phi(x + h_2) + \phi(x)\]
for any $x$ such that $x, x+ h_1, x+ h_2, x+ h_1 + h_2 \in B_N$.

It turns out that for the purposes of analysing Type I and II sums the cutoff $1_{B_N}$ plays a subservient r\^ole. The phase $\phi$, on the other hand, is crucial. The bulk of \cite{green-tao-u3mobius} is devoted to showing that if either a Type I or a Type II sum involving some $f$ as in \eqref{local-quad-def} is large, then $\phi$ is \emph{major arc}. This is a direct analogue of the proof of Davenport's bound as phrased at the end of \S \ref{vaughan-sec} (the ``inverse'' approach). Roughly speaking, $\phi$ is said to be major arc if $q \phi''(h_1,h_2)$ is small for some smallish $q$ and all $h_1,h_2$, which in turn essentially means that $\phi$ is slowly varying on $B_N$ intersected with any fixed progression $a \pmod{q}$. For a detailed discussion see \cite{green-tao-u3mobius}. Suffice it to say that the passage from large Type I/II sum to $\phi$ being major arc is long and difficult, and requires many applications of the Cauchy-Schwarz inequality to manipulate the phase $\phi$ into a helpful form, as well as basic tools of equidistribution such as a version of the Erd\H{o}s-Tur\'an inequality.

Recalling Proposition \ref{prop2}, one has reduced the case $s = 2$ of Conjecture \ref{mob-conj} to the statement that 
\[ \mathbb{E}_{n\leq N} \mu(n) 1_{B_N}(n) e(\phi(n)) \ll_A \log^{-A} N\]
for any major arc phase $\phi$. It turns out that $1_{B_N}(n) e(\phi(n))$ can, in this case, be closely approximated by a sum of linear phases $e(\theta n)$, and so we may conclude using Proposition \ref{davenport-prop}.

Note that this analysis has the flavour of an induction on $s$, the step of the nilsequence we are considering. We expect to see this more clearly when addressing the general case of Conjecture \ref{mob-conj} in future work.

\section{Future directions}

The most obvious avenue of research left open is to generalise everything we have done for $s = 2$ to the case $s \geq 3$. In particular we would like inverse theorems for the $U^k$-norms for $k \geq 4$, and a proof of Conjecture \ref{mob-conj} for $s \geq 3$. We are currently working towards this goal. We expect that the methods of Gowers \cite{gowers-longaps} can be adapted to achieve the inverse theorem, though this will not be straightforward. It is also very likely that the ``inverse'' approach to handling Type I and II sums can be adapted to the higher-step case of Conjecture \ref{mob-conj}, though again we do not expect this to be wholly straightforward.

It would be very desirable to have good bounds for error terms such as the $o(1)$ in Theorem \ref{4-aps-asymptotic}. We are sure that our current estimate for the error in Theorem \ref{4-aps-asymptotic} is the worst that has ever featured in analytic number theory -- the error term is a \emph{completely ineffective} $o(1)$! Ultimately this is because to show that the error is less than $\delta$ one finds oneself needing to rule out a real zero of some $L(s,\chi)$, $\chi$ a primitive quadratic character to the modulus $q$, with $s > 1 - C q^{-\epsilon}$, where $\epsilon = \epsilon(\delta) \rightarrow 0$ as $\delta \rightarrow 0$. Siegel's theorem states that for any $\epsilon > 0$ there is such a $C$, but it is, of course, not possible to specify $C$ effectively.

It is clear that the spectre of ineffectivity does not rear its head under the assumption of GRH, and we believe that our methods lead to an error term of the form $\log^{-c} N$ in Theorem \ref{4-aps-asymptotic}. 

There are other, presumably more tractible, ways in which one might obtain an explicit error term. Improvements to the combinatorial tools used in \cite{green-tao-u3inverse}, particularly advances on the circle of conjectures known as the ``polynomial Freiman-Ruzsa conjecture'', could be very helpful here.

We turn now to goals which lie further away. I have hinted at various places in this survey that the way in which we see nilsequences arising is very long-winded and, presumably, not the ``right'' way. The ergodic theorists \cite{host-kra,ziegler}  do admittedly discover the r\^ole of these functions somewhat less painfully (albeit after setting up a good deal of notation). Nilsequences seem such natural objects, however, that there ought to be a much better way of appreciating their place in the study of systems of linear equations. Recalling that $\Vert f \Vert_{U^2}$ is essentially the $L^4$ norm of $\widehat{f}$ one might even ask, for example,

\begin{question}
Is there a usable ``formula'' relating $\Vert f \Vert_{U^3}$ and certain of the ``nil-fourier coefficients'' $\mathbb{E}_{n \leq N} f(n) F(T_g^n \cdot x)$?
\end{question}

Such a formula would assuredly have to be very exotic on account of the vast profusion of nilsequences which might enter into consideration. The nilsequences are not naturally parametrised by anything so simple as the circle $S^1$, which gave its name to the classical circle method.

Let us conclude with some speculations on non-linear systems of equations, where our knowledge is at present essentially non-existent. We have seen in Conjecture \ref{conj} that the behaviour of an any system $\mathbf{A} \mathbf{x} = \mathbf{b}$, where $\mathbf{A}$ is non-degenerate in the sense of Definition \ref{non-degenerate-def}, should be governed by a very ``hard'' or ``algebraic'' collection of \emph{characteristic functions,} in this case the nilsequences. 

On the other hand degenerate linear systems, such as $x_1 - x_2 = 1$, do not have this property. To see this, suppose that $N = 2m$ is even and let $A \subseteq [N]$ be a set formed by setting $A \cap \{2i, 2i+1\} = \{2i\}$ or $\{2i+1\}$, these choices being independent in $i$ for $i = 0,\dots, m-1$. Then $|A| = N/2$, and $A$ is indistinguishable from a truly random set by taking inner products with any conceivable ``hard'' character such as a linear or quadratic phase. However, $A$ is expected to have about $N/8$ solutions to $x_1 - x_2 = 1$, whereas a random set has about twice this many.

One might call an equation or system of equations for which a ``hard'' characteristic system exists a \emph{mixing} system. We do not have a precise definition of this notion. Some non-linear equations are known to be mixing -- for example, the linear phases $e(\theta n)$ form a characteristic system for the equation $x_1 + x_2 = x_3^2$. Many more are not. It would be very interesting to know, for example, whether the equation $x_1 x_2 - x_3 x_4 = 1$ is mixing and, if so, what a characteristic system for it might be. This seems to be a very difficult question as the analysis of this equation even in very specific situations involves deep methods from the theory of automorphic forms.


\frenchspacing
\begin{thebibliography}{99}


\bibitem{baker-harman} Baker, R.~C. and Harman, G., Exponential sums formed with the M\"obius function, \emph{J. London Math. Soc.} (2) \textbf{43} (1991), no. 2, 193--198.

\bibitem{balog} Balog, A., Linear equations in primes, \emph{Mathematika} \textbf{39} (1992), 367--378.

\bibitem{bourgain} Bourgain, J., On triples in arithmetic progression, \emph{GAFA} \textbf{9} (1999), 968--984.

\bibitem{chowla} Chowla, S., There exists an infinity of 3---combinations of primes in A.P., \emph{Proc. Lahore. Philos. Soc.} \textbf{6} (1944), no. 2, 15--16.

\bibitem{davenport-book} Davenport, H., Multiplicative number theory, Third edition. \emph{Graduate Texts in Mathematics} \textbf{74}. Springer-Verlag, New York, 2000. xiv+177 pp. 

\bibitem{furst-weiss}
Furstenberg, H. and Weiss, B., A mean ergodic theorem for $1/N \sum_{n=1}^N f(T^n x) g(T^{n^2} x)$, in \emph{Convergence in ergodic theory and probability (Columbus OH 1993)}, 193--227, Ohio State Univ. Math. Res. Inst. Publ., 5. de Gruyter, Berlin, 1996.

\bibitem{gy} Goldston, D.~A. and Y{\i}ld{\i}r{\i}m, C.~Y., Small gaps between primes, I, \emph{preprint}.

\bibitem{gowers-4aps} 
Gowers, W.~T., A new proof of Szemer\'edi's theorem for arithmetic 
progressions of length four, \emph{GAFA} \textbf{8} (1998), 529--551.

\bibitem{gowers-icm} Gowers, W.~T., Fourier analysis and Szemer\'edi's theorem, in \emph{Proceedings of the Inetrnational Congress of Mathematicians, Berlin 1998}, Vol. 1.
 
\bibitem{gowers-longaps}
Gowers, W.~T., A new proof of Szemer\'edi's theorem, \emph{GAFA} \textbf{11} (2001), 465-588.

\bibitem{green-survey} Green, B.~J., Long arithmetic progressions of primes, \emph{preprint,} submitted to Proceedings of the Gauss-Dirichlet Conference, G\"ottingen 2005.

\bibitem{green-tao-longprimeaps} Green, B.~J. and Tao, T.~C., The primes contain arbitrarily long arithmetic progressions, \emph{to appear, Annals of Mathematics.}

\bibitem{green-tao-u3inverse} Green, B.~J. and Tao, T.~C., An inverse theorem for the Gowers $U^3$-norm, with applications, \emph{submitted}.

\bibitem{green-tao-u3mobius} Green, B.~J. and Tao, T.~C., Quadratic uniformity of the M\"obius function, \emph{preprint}.

\bibitem{green-tao-prime4aps} Green, B.~J. and Tao, T.~C., Two linear equations in four prime unknowns, \emph{in preparation}.

\bibitem{hardy-littlewood-primes1} Hardy, G.~H. and Littlewood, J.~E., Some problems of ``Partitio Numerorum''. III. On the expression of a number as a sum of primes, \emph{Acta. Math.} \textbf{44} (1923), 1--70.

\bibitem{hardy-littlewood-primes2} Hardy, G.~H. and Littlewood, J.~E., Some problems of ``Partitio Numerorum''. V. A further contribution to the study of Goldbach's problem, \emph{Proc. London Math. Soc.} (2) \textbf{22} (1923), 46--56.

\bibitem{hardy-ramanujan} Hardy, G.~H. and Ramanujan, S., Asymptotic formul{\ae} in combinatory analysis, \emph{Proc. London Math. Soc.} (2) \textbf{17} (1918), 75--115. 

\bibitem{heath-brown}
Heath-Brown, D.~R. Three primes and an almost prime in arithmetic progression, \emph{J. London Math. Soc.} (2) \textbf{23} (1981), 396--414.

\bibitem{host-kra} Host, B. and Kra, B. Non-conventional ergodic averages and nilmanifolds, \emph{Annals of Mathematics} \textbf{161} (2005), no. 1, 397--488.

\bibitem{iwaniec-kowalski} Iwaniec, H. and Kowalski, E. Analytic number theory, \emph{AMS Colloq. Publ.} \textbf{53}, AMS, Providence 2004.

\bibitem{iwaniec-luo-sarnak} Iwaniec, H., Luo, W and Sarnak, P., Low lying zeroes of families of $L$-functions, \emph{IHES Publ. Math.} \textbf{91} (2000), 55--131.

\bibitem{kra-survey} Kra, B., The Green-Tao Theorem on arithmetic progressions in the primes: an ergodic point of view, \emph{Bull. Amer. Math. Soc.} \textbf{43} (2006), 3--23.

\bibitem{kra-icm} Kra, B., From combinatorics to ergodic theory and back again, \emph{Proceedings of ICM 2006, Madrid}.

\bibitem{kumchev-survey} Kumchev, A.~V. and Tolev, D.~I., An invitation to additive prime number theory,  \emph{Serdica Math. J.}  \textbf{31}  (2005),  no. 1-2, 1--74.

\bibitem{salem-zygmund} Salem, R. and Zygmund, A., Some properties of trigonometric series whose terms have random signs, \emph{Acta Math.} \textbf{91} (1954), 245--301.

\bibitem{tao-survey-1} Tao, T.~C., Arithmetic progressions and the primes -- El Escorial Lectures, \emph{to appear, 2004 El Escorial proceedings}.

\bibitem{tao-survey-2} Tao, T.~C., Obstructions to uniformity, and arithmetic patterns in the primes, \emph{preprint}.

\bibitem{vdC} Van der Corput, J.~G., \"Uber Summen von Primzahlen und Primzahlquadraten, \emph{Math. Ann.} \textbf{116} (1939), 1--50.

\bibitem{vaughan-cr} Vaughan, R.~C., Sommes trigonom\'etriques sur les nombres premiers, \emph{C. R. Acad. Sci. Paris S\'er. A-B} \textbf{285} (1977), no. 16, A981--A983.

\bibitem{vaughan-legacy} Vaughan, R.~C., Hardy's Legacy to Number Theory, \emph{J. Austral. Math. Soc. (Series A)} \textbf{65} (1998), 238--266.

\bibitem{waring-survey} Vaughan, R.~C. and Wooley, T.D., Waring's problem: a survey, in \emph{Number theory for the millennium, III (Urbana, IL, 2000),}  301--340, A K Peters, Natick, MA, 2002. 

\bibitem{vinogradov-paper} Vinogradov, I.~M., Representation of an odd number as the sum of three primes, \emph{Dokl. Akad. Nauk SSSR} \textbf{15} (1937), 291--294.

\bibitem{wooley-survey} Wooley, T.~D., Diophantine problems in many variables: the r\^ole of additive number theory, in \emph{Topics in Number Theory, S. D. Ahlgren et al. (eds.)}, Kluwer Academic Publishers, 1999, pp. 49-83.

\bibitem{ziegler} Ziegler, T., Universal characteristic factors and Furstenberg averages, \emph{to appear, J. Amer. Math. Soc.}

 

\end{thebibliography}
\end{document}